# MORE ON HALFWAY NEW CARDINAL CHARACTERISTICS

BARNABÁS FARKAS, LUKAS DANIEL KLAUSNER, AND MARC LISCHKA

Abstract. We continue investigating variants of the splitting and reaping numbers introduced in [BHK+23]. In particular, answering a question raised there, we prove the consistency of $\operatorname{cof}(\mathcal{M}) < \mathfrak{s}_{1/2}$ and of $\mathfrak{r}_{1/2} < \operatorname{add}(\mathcal{M})$. Moreover, we discuss their natural generalisations $\mathfrak{s}_\rho$ and $\mathfrak{r}_\rho$ for $\rho \in (0,1)$, and show that $\mathfrak{r}_\rho$ does not depend on $\rho$.

## 1. Introduction

Let us recall the classical splitting number: Given $S, R \in [\omega]^\omega$ we say that $S$ *splits* $R$, in notation $S \mid R$, if $|S \cap R| = |R \smallsetminus S| = \omega$; and the *splitting number* is defined as

$$\mathfrak{s} = \min\{|\mathcal{S}| \mid \mathcal{S} \subseteq [\omega]^\omega \text{ and } \forall\, R \in [\omega]^\omega\ \exists\, S \in \mathcal{S}\colon S \mid R\}.$$

In [BHK+23], among many other new cardinal invariants, the following variant of $\mathfrak{s}$ was introduced: For $S, R \in [\omega]^\omega$ we say that $S$ *bisects* $R$, written as $S \mid_{1/2} R$, if

$$\frac{|S \cap R \cap n|}{|R \cap n|} \xrightarrow{n \to \infty} \frac{1}{2},$$

and $\mathfrak{s}_{1/2}$ is defined by replacing $\mid$ with $\mid_{1/2}$ in the definition of $\mathfrak{s}$,

$$\mathfrak{s}_{1/2} = \min\{|\mathcal{S}| \mid \mathcal{S} \subseteq [\omega]^\omega \text{ and } \forall\, R \in [\omega]^\omega\ \exists\, S \in \mathcal{S}\colon S \mid_{1/2} R\}.$$

As $S \mid_{1/2} R$ implies $S \mid R$, $\mathfrak{s} \leq \mathfrak{s}_{1/2}$ immediately follows. It also turned out that $\operatorname{cov}(\mathcal{M}) \leq \mathfrak{s}_{1/2} \leq \operatorname{non}(\mathcal{N})$, where $\operatorname{cov}(\mathcal{M})$ stands for the covering number of the meagre ideal and $\operatorname{non}(\mathcal{N})$ for the uniformity of the null ideal (see below). Moreover, two of these inequalities are consistently strict (see [BHK+23, end of Section 2]): $\mathfrak{s} < \mathfrak{s}_{1/2}$ in the Cohen model and $\operatorname{cov}(\mathcal{M}) < \mathfrak{s}_{1/2}$ in the Mathias model. One of the most interesting remaining open questions is the following:

**Question A.** *Is* $\mathfrak{s}_{1/2} < \operatorname{non}(\mathcal{N})$ *consistent?*

Of course, similar questions were raised regarding separations of $\mathfrak{s}_{1/2}$ and other classical invariants. For example (see [BHK+23, Question A (Q2)]), motivated by the fact that $\mathfrak{s} \leq \mathfrak{d}$, where $\mathfrak{d}$ stands for the dominating number (see below), it is natural to ask the following:

**Question B.** *Is* $\mathfrak{d} < \mathfrak{s}_{1/2}$ *consistent?*

In section 2, we present a short overview on the relevant cardinal invariants (from [BHK+23]) and recall the inequalities and consistently strict inequalities between variants of $\mathfrak{s}$ and other classical cardinal characteristics.

In section 3, by dualising the results from section 2, we outline the inequalities between the variants of the reaping number $\mathfrak{r}$ and show the consistency of almost all possible strict inequalities in the dual diagram.

2020 *Mathematics Subject Classification.* Primary 03E17; Secondary 03E40.
The first author was supported by the Austrian Science Fund (FWF) projects P29907 "Borel Ideals and Filters" and I5918 "Analytic $P$-Ideals, Banach Spaces, and Measure Algebras".

In section 4, answering Question B positively, we present two models of $\operatorname{cof}(\mathcal{M}) < \mathfrak{s}_{1/2}$: (1) the $\omega_2$-stage countable support iteration of a modified infinitely equal forcing and (2) the dual Hechler model, that is, the $\omega_1$-stage finite support iteration of the Hechler forcing over a model of $\textsc{ma} + \mathfrak{c} \geq \omega_2$.

In section 5, we define natural generalisations of $\mathfrak{s}_{1/2}$ and $\mathfrak{r}_{1/2}$, namely $\mathfrak{s}_\rho$ and $\mathfrak{r}_\rho$ for $\rho \in (0,1)$. We discuss their lower and upper bounds and show that $\mathfrak{r}_\rho$ does not depend on $\rho$.

## 2. Variants of $\mathfrak{s}$ and $\mathfrak{r}$

When studying cardinal characteristics the following framework can come handy, especially when dualising inequalities and in the context of forcing. For a general overview on this framework see e. g. [Bla10] or [Fre15, Section 512]. (Here, we mostly follow the notation of the second reference.)

**Definition 2.1.** A *relational system* is a triplet $\mathbf{R} = (X, \sqsubset, Y)$ where (i) $\sqsubset$ is a relation on $X \times Y$, (ii) $\operatorname{dom}(\sqsubset) = X$ and (iii) there is no single $y \in Y$ such that $x \sqsubset y$ for every $x \in X$.

A set $U \subseteq X$ is $\mathbf{R}$-*unbounded* if there is no single $y \in Y$ $\sqsubset$-above all elements of $U$; a set $D \subseteq Y$ is $\mathbf{R}$-*dominating* if for every $x \in X$, there is a $y \in D$ $\sqsubset$-above $x$. The *(un)bounding* and *dominating numbers* of $\mathbf{R}$ are

$$\mathfrak{b}(\mathbf{R}) = \mathfrak{b}(\sqsubset) = \min\{|U| \mid U \subseteq X \text{ is } \mathbf{R}\text{-unbounded}\},$$
$$\mathfrak{d}(\mathbf{R}) = \mathfrak{d}(\sqsubset) = \min\{|D| \mid D \subseteq Y \text{ is } \mathbf{R}\text{-dominating}\}.$$

The *dual of* $\mathbf{R}$ is the relational system $\mathbf{R}^\perp = (Y, \not\sqsupset, X)$ (which satisfies conditions (i)–(iii) automatically). Clearly, $\mathfrak{b}(\mathbf{R}^\perp) = \mathfrak{d}(\mathbf{R})$ and $\mathfrak{d}(\mathbf{R}^\perp) = \mathfrak{b}(\mathbf{R})$.

**Example 2.2.** We recall all relational systems and cardinal characteristics we need in the first four sections:

(1) Let $\mathbf{Dom} = (\omega^\omega, \leq^*, \omega^\omega)$, where $f \leq^* g$ if $\{n \mid g(n) < f(n)\}$ is finite. Then $\mathfrak{b} = \mathfrak{b}(\mathbf{Dom})$ and $\mathfrak{d} = \mathfrak{d}(\mathbf{Dom})$.

(2) Let $\mathcal{I}$ be an ideal on an infinite set $X$, that is, $[X]^{<\omega} \subseteq \mathcal{I}$ and $C \in \mathcal{I}$ whenever $C \subseteq A \cup B$ for some $A, B \in \mathcal{I}$; we will always assume that $X \notin \mathcal{I}$. Consider the relational systems $\mathbf{Cof}(\mathcal{I}) = (\mathcal{I}, \subseteq, \mathcal{I})$ and $\mathbf{Cov}(\mathcal{I}) = (X, \in, \mathcal{I})$. The four classical cardinal invariants of $\mathcal{I}$ are the following cardinals:

$$\operatorname{add}(\mathcal{I}) = \mathfrak{b}(\mathbf{Cof}(\mathcal{I})), \qquad \operatorname{cof}(\mathcal{I}) = \mathfrak{d}(\mathbf{Cof}(\mathcal{I})),$$
$$\operatorname{non}(\mathcal{I}) = \mathfrak{b}(\mathbf{Cov}(\mathcal{I})), \qquad \operatorname{cov}(\mathcal{I}) = \mathfrak{d}(\mathbf{Cov}(\mathcal{I})).$$

We know that $\operatorname{add}(\mathcal{I}) \leq \operatorname{cov}(\mathcal{I}), \operatorname{non}(\mathcal{I}) \leq \operatorname{cof}(\mathcal{I})$ always holds. We are particularly interested in two specific ideals, namely $\mathcal{M} = \{\text{meagre subsets of } 2^\omega\}$ and $\mathcal{N} = \{\text{null subsets of } 2^\omega\}$. The well-known Cichoń's diagram can be summarised in the following Tukey connections (see below): $\mathbf{Cov}(\mathcal{M}) \preceq \mathbf{Dom} \preceq \mathbf{Cof}(\mathcal{M}) \preceq \mathbf{Cof}(\mathcal{N})$ and $\mathbf{Cov}(\mathcal{N}) \preceq \mathbf{Cov}(\mathcal{M})^\perp$ (and the facts that $\operatorname{add}(\mathcal{M}) = \min\{\mathfrak{b}, \operatorname{cov}(\mathcal{M})\}$ and $\operatorname{cof}(\mathcal{M}) = \max\{\operatorname{non}(\mathcal{M}), \mathfrak{d}\}$).

(3) Let $\mathbf{Reap} = ([\omega]^\omega, \nmid, [\omega]^\omega)$ where $S \mid X$ ("$S$ splits $X$") if $|S \cap X| = |X \setminus S| = \omega$. Then $\mathfrak{s} = \mathfrak{b}(\mathbf{Reap})$ and $\mathfrak{r} = \mathfrak{d}(\mathbf{Reap})$.

(4) Let $\mathbf{Reap}_{1/2} = ([\omega]^\omega, \nmid_{1/2}, [\omega]^\omega)$, where $S \mid_{1/2} X$ ("$S$ bisects $X$") if

$$\frac{|S \cap X \cap n|}{|X \cap n|} \xrightarrow{n \to \infty} \frac{1}{2}.$$

Then $\mathfrak{s}_{1/2} = \mathfrak{b}(\mathbf{Reap}_{1/2})$ and $\mathfrak{r}_{1/2} = \mathfrak{d}(\mathbf{Reap}_{1/2})$.

(5) For $\varepsilon \in (0, {}^1\!/\!_2)$ let $\mathbf{Reap}_{1/2\pm\varepsilon} = ([\omega]^\omega, \not|_{1/2\pm\varepsilon}, [\omega]^\omega)$, where $S \mid_{1/2\pm\varepsilon} X$ ("$S$ $\varepsilon$-almost bisects $X$") if for all but finitely many $n$

$$\frac{1}{2} - \varepsilon < \frac{|S \cap X \cap n|}{|X \cap n|} < \frac{1}{2} + \varepsilon.$$

Then $\mathfrak{s}_{1/2\pm\varepsilon} = \mathfrak{b}(\mathbf{Reap}_{1/2\pm\varepsilon})$ and $\mathfrak{r}_{1/2\pm\varepsilon} = \mathfrak{d}(\mathbf{Reap}_{1/2\pm\varepsilon})$.

Unfortunately, it is still unclear whether $\mathfrak{s}_{1/2\pm\varepsilon}$ depends on $\varepsilon$ (see [BHK+23, Question A (Q3)]). Whenever we claim anything about $\mathbf{Reap}_{1/2\pm\varepsilon}$ or its invariants, we mean that our claim holds for every $\varepsilon \in (0, {}^1\!/\!_2)$.

**Theorem 2.3** (see [BHK+23, Theorem 2.4])**.** *The following relations hold, where $a \longrightarrow b$ means "$a \leq b$, consistently $a < b$" and $a \dashrightarrow b$ means "$a \leq b$, possibly $a = b$":*

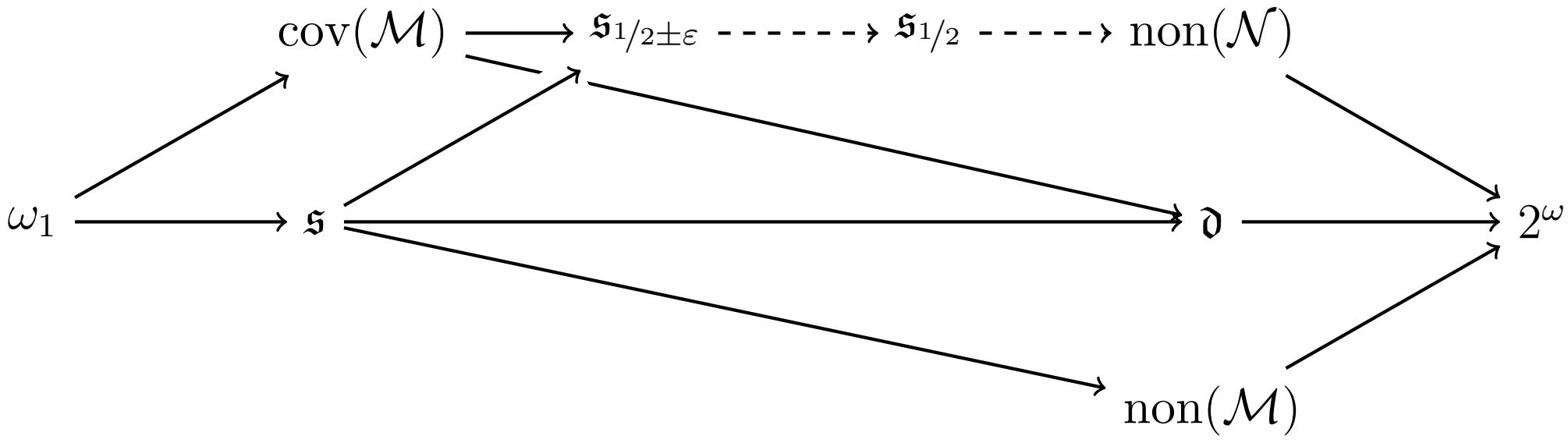

Regarding further inequalities between these cardinals (apart from separating $\mathfrak{s}_{1/2}$ from $\mathfrak{s}_{1/2\pm\varepsilon}$ and/or non($\mathcal{N}$)), there was only one question left open, namely, if $\mathfrak{s}_{1/2} \leq \mathfrak{d}$ or $\mathfrak{s}_{1/2\pm\varepsilon} \leq \mathfrak{d}$ hold. To show that this inequality does not hold, in section 4 we will present a construction based on countable support iteration as well as one based on finite support iteration.

## 3. The Dual Diagram

Although in the proof of Theorem 2.3 it was not stated explicitly, all inequalities were proved (or could have easily been proved) via Tukey connections.

**Definition 3.1.** Let $\mathbf{R}_0 = (X_0, \sqsubset_0, Y_0)$ and $\mathbf{R}_1 = (X_1, \sqsubset_1, Y_1)$ be relational systems. A pair $(F, G)$ of functions $F \colon X_0 \to X_1$ and $G \colon Y_1 \to Y_0$, in short $(F, G) \colon \mathbf{R}_0 \to \mathbf{R}_1$, is a *(Galois-)Tukey connection* from $\mathbf{R}_0$ to $\mathbf{R}_1$ if for any $x_0 \in X_0$ and $y_1 \in Y_1$, $F(x_0) \sqsubset_1 y_1$ implies $x_0 \sqsubset_0 G(y_1)$. This is commonly visualised as follows:

$$\begin{array}{ccccc} G(y_1) \in Y_0 & \xleftarrow{\;G\;} & Y_1 \ni & y_1 \\ \sqsubset_0 & \Longleftarrow & & \sqsubset_1 \\ x_0 \;\in X_0 & \xrightarrow[F]{} & X_1 \ni & F(x_0) \end{array}$$

If there is a Tukey connection from $\mathbf{R}_0$ to $\mathbf{R}_1$, we say that $\mathbf{R}_0$ is *Tukey-below* $\mathbf{R}_1$ and write $\mathbf{R}_0 \preccurlyeq \mathbf{R}_1$. Note that $(F, G) \colon \mathbf{R}_0 \to \mathbf{R}_1$ is a Tukey connection iff $(G, F) \colon \mathbf{R}_1^\perp \to \mathbf{R}_0^\perp$ is a Tukey connection. We say that $\mathbf{R}_0$ and $\mathbf{R}_1$ are *Tukey-equivalent*, denoted by $\mathbf{R}_0 \equiv \mathbf{R}_1$, when $\mathbf{R}_0 \preccurlyeq \mathbf{R}_1$ and $\mathbf{R}_1 \preccurlyeq \mathbf{R}_0$ both hold.

Recall that $\mathbf{R}_0 \preccurlyeq \mathbf{R}_1$ implies that $\mathfrak{b}(\mathbf{R}_0) \geq \mathfrak{b}(\mathbf{R}_1)$ and $\mathfrak{d}(\mathbf{R}_0) \leq \mathfrak{d}(\mathbf{R}_1)$. For example, the proof of $\mathfrak{s}_{1/2} \leq$ non($\mathcal{N}$) in Theorem 2.3 actually shows $\mathbf{Cov}(\mathcal{N}) \preceq \mathbf{Reap}_{1/2}$ etc., and hence the dual inequalities (e. g. cov($\mathcal{N}$) $\leq \mathfrak{r}_{1/2}$) hold as well.

Let us point out that, after appropriate coding, all relational systems and Tukey connections we discuss in this paper are Borel in the following sense: $\mathbf{R} = (X, \sqsubset, Y)$ is Borel if $X$ and $Y$ are Polish spaces and $\sqsubset$ is a Borel subset of $X \times Y$; a Tukey reduction $(F, G)\colon (X_0, \sqsubset_0, Y_0) \to (X_1 \sqsubset_1, Y_1)$ between Borel systems is Borel if both $F$ and $G$ are Borel functions. This can be particularly useful in the context of forcing: If $\mathbf{R}$ is Borel, then we say that a forcing notion $\mathbb{P}$ is $\mathbf{R}$-*dominating* if

$$V^{\mathbb{P}} \vDash \exists y \in Y \, \forall x \in X \cap V \colon x \sqsubset y.$$

It follows that if there is a Borel Tukey connection $\mathbf{R}_0 \to \mathbf{R}_1$ and $\mathbb{P}$ is $\mathbf{R}_1$-dominating, then $\mathbb{P}$ is $\mathbf{R}_0$-dominating as well. Instead of "$\mathbf{R}$-dominating", we will use the common more specific terms, e. g. "$\mathbb{P}$ adds a dominating real (over $V$)" means that $\mathbb{P}$ is $\mathbf{Dom}$-dominating, "$\mathbb{P}$ adds a random real" means that $\mathbb{P}$ is $\mathbf{Cov}(\mathcal{N})^{\perp}$-dominating, etc.; we can talk about adding splitting, bisecting or $\varepsilon$-almost bisecting reals (over $V$) analogously.

To dualise Theorem 2.3, we hence have to check the consistency of the strict inequalities. Before stating the dual form of Theorem 2.3, let us take a closer look at $\mathfrak{r}_{1/2}$. The Borel reducibility $\mathbf{Reap}_{1/2}^{\perp} \preceq \mathbf{Cov}(\mathcal{N})^{\perp}$ (because $\mathbf{Cov}(\mathcal{N}) \preceq \mathbf{Reap}_{1/2}$, see above) implies that if we add random reals, we also add bisecting reals over $V$. The next two lemmas illustrate that "too tame" forcing notions cannot increase $\mathfrak{r}_{1/2}$. We assume that all forcing notions are atomless and separative; in particular, we assume that every condition $p$ in a forcing notion $\mathbb{P}$ has incompatible extensions.

**Lemma 3.2.** *If $\mathbb{P}$ is $\sigma$-centred, then it cannot add $\varepsilon$-almost bisecting reals.*

*Proof.* Let $\mathbb{P} = \bigcup_{n \in \omega} C_n$, where each $C_n$ is centred, let $\dot{B}$ be a $\mathbb{P}$-name for an element of $[\omega]^{\omega}$, and assume towards a contradiction that

$$p \Vdash \text{“}\dot{B} \text{ } \varepsilon\text{-almost bisects every } X \in [\omega]^{\omega} \cap V\text{”}$$

for some $p \in \mathbb{P}$ and $\varepsilon \in (0, 1/2)$, that is,

$$p \Vdash \text{“}\frac{1}{2} - \varepsilon < \frac{|\dot{B} \cap X \cap n|}{|X \cap n|} < \frac{1}{2} + \varepsilon \text{ for almost all } n\text{”}$$

for every $X \in [\omega]^{\omega} \cap V$. Fix an interval partition $(I_n)$ of $\omega$ such that $|I_0| \geq 2$ and $|I_n| > 2^n |I_{<n}|$ for every $n$ (where $I_{<n} = \bigcup_{k<n} I_k$ and $I_{<0} = \varnothing$). First of all, we show that

$$p \Vdash \text{“}\frac{1}{2} - \varepsilon' < \frac{|\dot{B} \cap I_n|}{|I_n|} < \frac{1}{2} + \varepsilon' \text{ for almost all } n\text{”}$$

holds for every $\varepsilon' \in (\varepsilon, 1/2)$. To see this, fix such an $\varepsilon'$; then, as $p \Vdash$"$\dot{B}$ $\varepsilon$-almost bisects $\omega$", $p$ forces that for every sufficiently large $n$,

$$\begin{aligned}\frac{|\dot{B} \cap I_n|}{|I_n|} &\leq \frac{(1 + 2^{-n})|\dot{B} \cap I_{\leq n}|}{(1 + 2^{-n})|I_n|} < \frac{(1 + 2^{-n})|\dot{B} \cap I_{\leq n}|}{|I_{\leq n}|} \\ &< (1 + 2^{-n}) \left( \frac{1}{2} + \varepsilon \right) < \frac{1}{2} + \varepsilon',\end{aligned}$$

where the second inequality follows from $|I_n| + 2^{-n}|I_n| > |I_n| + |I_{<n}|$. The lower bound can be established similarly: $p$ forces that for every sufficiently large $n$,

$$\frac{|\dot{B} \cap I_n|}{|I_n|} > \frac{|\dot{B} \cap I_n|}{|I_n|} + \frac{|I_{<n}|}{|I_n|} - 2^{-n} \geq \frac{|\dot{B} \cap I_{\leq n}|}{|I_{\leq n}|} - 2^{-n} > \frac{1}{2} - \varepsilon - 2^{-n} > \frac{1}{2} - \varepsilon'.$$

Fix an $\varepsilon'$ as above. Then there is a $p' \leq p$ and an $N \in \omega$ such that

$$p' \Vdash \text{“}\frac{1}{2} - \varepsilon' < \frac{|\dot{B} \cap I_n|}{|I_n|} < \frac{1}{2} + \varepsilon' \text{ for every } n \geq N\text{”}.$$

By modifying $\dot{B}$ on $I_{<N}$ (which does not affect the property of being $\varepsilon$-almost bisecting), we can assume that $p'$ forces these inequalites for every $n$; in particular, $\dot{B} \cap I_n \neq \varnothing$ for every $n$ (as $\varepsilon'$ can be very small, we assume that $|I_n|$ is even for every $n$). In the rest of the proof we assume that $p' = 1_{\mathbb{P}}$, that is,

$$(*_1) \qquad \Vdash \text{“}\frac{1}{2} - \varepsilon' < \frac{|\dot{B} \cap I_n|}{|I_n|} < \frac{1}{2} + \varepsilon' \text{ for every } n\text{”}.$$

Let $\mathcal{E}_n = \{E \subseteq I_n \mid \forall p \in C_n \; \exists q \leq p \colon q \Vdash \dot{B} \cap I_n = E\}$. Note that $\mathcal{E}_n \neq \varnothing$: Otherwise, for every $E \subseteq I_n$ we can fix $p_E \in C_n$ such that $p_E \Vdash \dot{B} \cap I_n \neq E$, but these $p_E$ have a common extension, which is a contradiction. Fix an arbitrary sequence $(E_n)_{n\in\omega} \in V$ such that $E_n \in \mathcal{E}_n$ for every $n$ (we know that $E_n \neq \varnothing$), and let $X = \bigcup_{n\in\omega} E_n \in [\omega]^\omega \cap V$.

Since $\dot{B}$ $\varepsilon$-almost bisects $X$, we can fix some $p \in \mathbb{P}$ and $M \in \omega$ such that

$$(*_2) \qquad p \Vdash \text{“}\frac{|\dot{B} \cap X \cap m|}{|X \cap m|} < \frac{1}{2} + \varepsilon \text{ for all } m \geq M\text{”}.$$

Now fix a $k$ such that

$$\min(I_{k+1}) \geq M \qquad \text{and} \qquad \frac{1}{\frac{2^{-k}}{1/2 - \varepsilon'} + 1} \geq \frac{1}{2} + \varepsilon.$$

As each condition has incompatible extensions, there are extensions of $p$ in infinitely many $C_n$; hence for some $n \geq k$, we can fix a $p' \in C_n$ below $p$. By the definition of $\mathcal{E}_n$, there is a $q \leq p'$ which forces that $\dot{B} \cap I_n = E_n = X \cap I_n$; in particular, $q$ forces that

$$\frac{|\dot{B} \cap X \cap I_{\leq n}|}{|X \cap I_{\leq n}|} \geq \frac{|\dot{B} \cap X \cap I_n|}{|I_{<n}| + |X \cap I_n|} = \frac{|E_n|}{|I_{<n}| + |E_n|} = \frac{1}{\frac{|I_{<n}|}{|E_n|} + 1}$$
$$> \frac{1}{\frac{|I_{<n}|}{(1/2 - \varepsilon')|I_n|} + 1} > \frac{1}{\frac{2^{-n}}{1/2 - \varepsilon'} + 1} \geq \frac{1}{2} + \varepsilon,$$

where the second inequality follows from Eq. ($*_1$) and the third inequality follows from $|I_{<n}| < 2^{-n}|I_n|$. This contradicts Eq. ($*_2$), because $I_{\leq n} = \min(I_{n+1}) \geq \min(I_{k+1}) \geq M$. ☐

Before the next lemma, let us recall the Laver property. A forcing notion $\mathbb{P}$ has the *Laver property* (see [BJ95, Definition 6.3.27]) if for every sequence $(H_n)_{n\in\omega}$ of non-empty finite sets and every $\mathbb{P}$-name $\dot{f} \in \prod_{n\in\omega} H_n$, there is an $S \in \prod_{n\in\omega} [H_n]^{2^n}$ (in $V$) such that $\Vdash_{\mathbb{P}}$“$\dot{f}(n) \in S(n)$ for almost all $n$”.

**Lemma 3.3.** *If $\mathbb{P}$ has the Laver property, then it cannot add $\varepsilon$-almost bisecting reals.*

*Proof.* Assume towards a contradiction that a $p \in \mathbb{P}$ forces that $\dot{B} \in V^{\mathbb{P}}$ $\varepsilon$-almost bisects every $X \in [\omega]^\omega \cap V$. As in the proof of Lemma 3.2, we fix an interval partition $(I_n)$ in $V$ such that $|I_0| \geq 2$ and $|I_n| > 2^n |I_{<n}|$ for every $n$ as well as an $\varepsilon' \in (\varepsilon, 1/2)$, and just like above, we can assume that

$$(*_3) \qquad p \Vdash \text{“}\frac{1}{2} - \varepsilon' < \frac{|\dot{B} \cap I_n|}{|I_n|} < \frac{1}{2} + \varepsilon' \text{ for every } n\text{”}.$$

Let $Q_0 = I_1$, $Q_1 = I_2 \cup I_3$, ..., $Q_m = \bigcup_{m'<2^m} I_{2^m+m'}$ be the union of the next $2^m$ many intervals in the partition. Applying the Laver property, there are an $S \in \prod_{m\in\omega} [\mathcal{P}(Q_m)]^{2^m}$ in $V$, $S(m) = \{S^m_{m'} \mid m' < 2^m\}$, a $p' \leq p$, and a $\mathbb{P}$-name $\dot{b}$ for an element of $\prod_{m\in\omega} 2^m$ such that $p' \Vdash$“$\dot{B} \cap Q_m = S^m_{\dot{b}(m)} \in S(m)$ for every $m$”. Define

$$X = \bigcup_{m\in\omega} \bigcup_{m'<2^m} S^m_{m'} \cap I_{2^m+m'} \in V.$$

Then $X$ is infinite because $X \cap I_{2^m+\dot{b}(m)} = S^m_{\dot{b}(m)} \cap I_{2^m+\dot{b}(m)} = \dot{B} \cap I_{2^m+\dot{b}(m)} \neq \varnothing$ for every $m$. We claim that $p' \Vdash$"$\dot{B}$ does not $\varepsilon$-almost bisect $X$". To see this, let $m$ be sufficiently large such that

$$\frac{1}{\frac{2^{-2^m}}{1/2-\varepsilon'}+1} \geq \frac{1}{2}+\varepsilon.$$

Then $p'$ forces that

$$\begin{aligned}\frac{|\dot{B}\cap X\cap I_{\leq 2^m+\dot{b}(m)}|}{|X\cap I_{\leq 2^m+\dot{b}(m)}|} &\geq \frac{|\dot{B}\cap I_{2^m+\dot{b}(m)}|}{|I_{<2^m+\dot{b}(m)}|+|\dot{B}\cap I_{2^m+\dot{b}(m)}|} = \frac{1}{\frac{|I_{<2^m+\dot{b}(m)}|}{|\dot{B}\cap I_{2^m+\dot{b}(m)}|}+1}\\ &> \frac{1}{\frac{|I_{<2^m+\dot{b}(m)}|}{(1/2-\varepsilon')|I_{2^m+\dot{b}(m)}|}+1} > \frac{1}{\frac{2^{-2^m-\dot{b}(m)}}{1/2-\varepsilon'}+1} \geq \frac{1}{2}+\varepsilon\end{aligned}$$

where we used Eq. ($*_3$) in the second inequality and $|I_{<n}|/|I_n| < 2^{-n}$ in the third one. $\square$

**Theorem 3.4.** *The following relations hold, where $a \longrightarrow b$ means "$a \leq b$, consistently $a < b$" and $a \dashrightarrow b$ means "$a \leq b$, possibly $a = b$", and there are no further provable inequalities between these cardinals:*

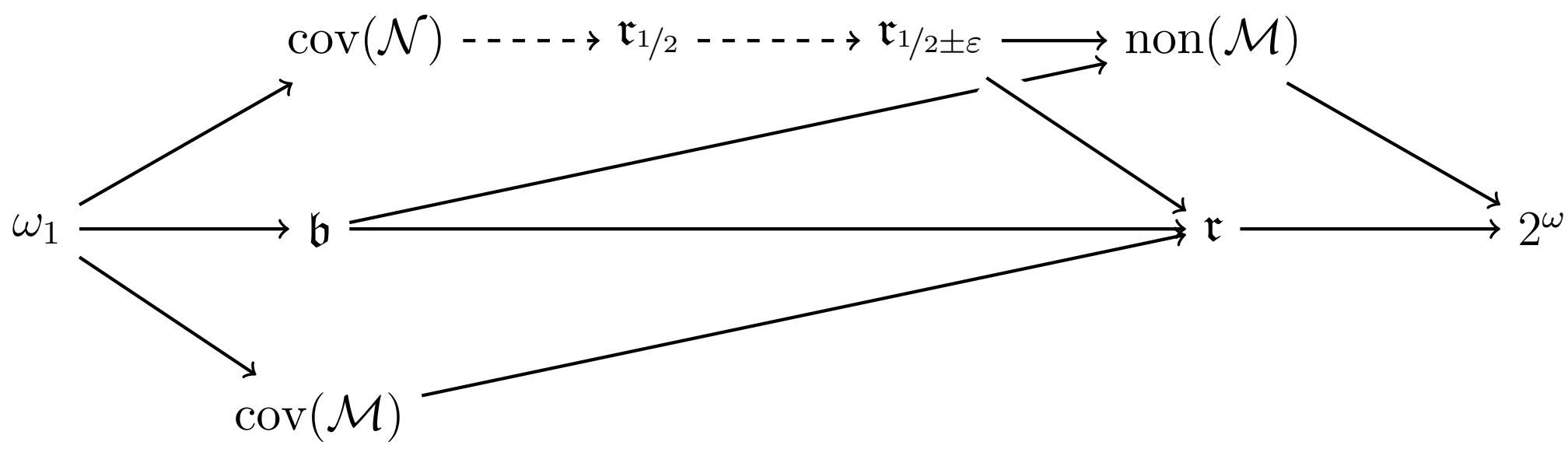

*Proof.* The inequalities follow by dualising the ones in Theorem 2.3. To show the consistency of the strict inequalities and that there are no further inequalities in ZFC, keeping in mind the consistent cuts of the Cichoń's diagram, it is enough to show that the following are consistent: $\mathfrak{r} < \mathrm{non}(\mathcal{M})$ and $\mathfrak{r}_{1/2\pm\varepsilon} < \mathfrak{b}$ (separately, of course).

$\mathfrak{r} < \mathrm{non}(\mathcal{M})$ holds in the model presented in [BS87, Section 5], because it is a model of $\mathfrak{u} < \mathfrak{s}$ and we know that $\mathfrak{r} \leq \mathfrak{u}$ and $\mathfrak{s} \leq \mathrm{non}(\mathcal{M})$ (for more details on the ultrafilter number $\mathfrak{u}$, see [Bla10, Section 9]).

We show that $\mathfrak{r}_{1/2\pm\varepsilon} < \mathrm{add}(\mathcal{M}) \leq \mathfrak{b}$ holds after the $\mathfrak{c}^+$ stage finite support iteration $\mathbb{D}_{\mathfrak{c}^+}$ of the Hechler forcing $\mathbb{D}$. We know that if $\kappa \geq \mathfrak{c}$ is regular, then $V^{\mathbb{D}_\kappa} \vDash \mathrm{add}(\mathcal{M}) = \min\{\mathfrak{b}, \mathrm{cov}(\mathcal{M})\} = \kappa = \mathfrak{c}$ because $\mathbb{D}$ adds both dominating and Cohen reals. Recall (folklore, see [GBG11]) that if $\delta \leq \mathfrak{c}$ (hence also if $\delta < \mathfrak{c}^+$) and $(\mathbb{P}_\alpha, \dot{\mathbb{Q}}_\beta)_{\alpha\leq\delta,\beta<\delta}$ is a finite support iteration of $\sigma$-centred forcing notions (that is, $\Vdash_\alpha$ "$\dot{\mathbb{Q}}_\alpha$ is $\sigma$-centred" for every $\alpha < \delta$), then $\mathbb{P}_\delta$ is also $\sigma$-centred. Applying Lemma 3.2, it follows that if $\delta \leq \mathfrak{c}^+$ and $(\mathbb{P}_\alpha, \dot{\mathbb{Q}}_\beta)_{\alpha\leq\delta,\beta<\delta}$ is a finite support iteration of $\sigma$-centred forcing notions, then

$$\Vdash_\delta \text{"no } S \in [\omega]^\omega \text{ can } \varepsilon\text{-almost bisect all elements of } [\omega]^\omega \cap V\text{";}$$

in particular, $V^{\mathbb{D}_{\mathfrak{c}^+}} \vDash \mathfrak{r}_{1/2\pm\varepsilon} \leq \mathfrak{c}^V < \mathrm{add}(\mathcal{M}) = \mathfrak{c}$. $\square$

Probably the most interesting remaining open question is the following:

**Question C.** *Is* $\mathrm{cov}(\mathcal{N}) < \mathfrak{r}_{1/2}$ *consistent?*

## 4. Models of $\mathrm{cof}(\mathcal{M}) < \mathfrak{s}_{1/2}$

Let us recall the infinitely equal forcing $\mathbb{EE}$ (see [BJ95, Definition 7.4.11]): $p \in \mathbb{EE}$ if $p$ is a function, $\mathrm{dom}(p) \subseteq \omega$ is coinfinite, and $p(n)\colon n \to 2$ is a function for every $n$; if $p$ and $q$ are conditions, $q \leq p$ if $q \supseteq p$. We know (see [BJ95, Lemma 7.4.12 and Lemma 7.4.14]) that $\mathbb{EE}$ and its countable support iterations are proper and $\omega^\omega$-bounding and preserve non-meagre sets (see also [BJ95, Lemma 6.3.21, Theorem 6.1.13 and Theorem 6.3.5]).

We will work with a variant of $\mathbb{EE}$. Basically, we do the following: (i) We allow longer characteristic functions in the $n$-th coordinate (but still below a common bound); and to adapt the forcing notion to our needs, (ii) we switch from characteristic functions to subsets and (iii) by shifting the underlying set in the $n$-th coordinate, we make sure that in the generic sequence, these finite sets are ordered consecutively.

More precisely, we fix an interval partition $(I_n)$ of $\omega$ as in Lemma 3.2 and Lemma 3.3, that is, $|I_0| \geq 2$ and $|I_n| > 2^n |I_{<n}|$ for every $n$, and define $\mathbb{P} = \mathbb{P}_{(I_n)}$ as follows: $p \in \mathbb{P}$ if $p$ is function such that

(a) $\mathrm{dom}(p) \subseteq \omega$ is coinfinite and
(b) $p(n) \subseteq I_n$ for every $n \in \mathrm{dom}(p)$;

$q \leq p$ if $q \supseteq p$. The very same proofs that work for $\mathbb{EE}$ show that $\mathbb{P}$ and its countable support iterations are proper and $\omega^\omega$-bounding and preserve non-meagre sets. It follows that if CH holds in the ground model, then $\mathbb{P}_{\omega_2}$ (the $\omega_2$-stage countable support iteration of $\mathbb{P}$) forces that $\mathrm{cof}(\mathcal{M}) = \max\{\mathrm{non}(\mathcal{M}), \mathfrak{d}\} = \omega_1$.

**Lemma 4.1.** *Let $\dot{G}$ be the canonical $\mathbb{P}$-name of the $\mathbb{P}$-generic filter and $\dot{X}$ be a $\mathbb{P}$-name such that*

$$\Vdash_{\mathbb{P}} \dot{X} = \bigcup \mathrm{ran}\Big(\bigcup \dot{G}\Big) = \bigcup \Big\{ p(n) \;\Big|\; p \in \dot{G} \text{ and } n \in \mathrm{dom}(p) \Big\} \; (\in [\omega]^\omega).$$

*Then no $S \in [\omega]^\omega \cap V$ can $\varepsilon$-almost bisect $\dot{X}$.*

*Proof.* Fix an $S \in [\omega]^\omega$, an $\varepsilon \in (0, 1/2)$ and a $p \in \mathbb{P}$. Pick an $n$ in $\omega \setminus \mathrm{dom}(p)$ such that

$$2^{-n+1} \leq \frac{1}{2} - \varepsilon \qquad \text{and hence} \qquad \frac{1}{2^{-n+1}+1} \geq \frac{1}{2} + \varepsilon.$$

We distinguish two cases:

<u>Case 1</u>: $|S \cap I_n| > |I_n|/2$. Define $q_n \in \mathbb{P}$, $\mathrm{dom}(q_n) = \mathrm{dom}(p) \cup \{n\}$, $q_n{\restriction}_{\mathrm{dom}(p)} = p$ and $q_n(n) = S \cap I_n$. Then $q_n \leq p$ and $q_n$ forces that

$$\begin{aligned}\frac{|S \cap \dot{X} \cap I_{\leq n}|}{|\dot{X} \cap I_{\leq n}|} &\geq \frac{|S \cap I_n|}{|I_{<n}| + |S \cap I_n|} = \frac{1}{\frac{|I_{<n}|}{|S \cap I_n|} + 1} \\ &> \frac{1}{\frac{|I_{<n}|}{|I_n|/2} + 1} > \frac{1}{2^{-n+1}+1} \geq \frac{1}{2} + \varepsilon.\end{aligned}$$

As $n$ can be arbitrarily large, $p$ cannot force that $S$ $\varepsilon$-almost bisects $\dot{X}$.

<u>Case 2</u>: $|S \cap I_n| \leq |I_n|/2$. Define $r_n \in \mathbb{P}$, $\mathrm{dom}(r_n) = \mathrm{dom}(p) \cup \{n\}$, $r_n{\restriction}_{\mathrm{dom}(p)} = p$ and $r_n(n) = I_n \setminus S$. Then $r_n \leq p$ and $r_n$ forces that

$$\frac{|S \cap \dot{X} \cap I_{\leq n}|}{|\dot{X} \cap I_{\leq n}|} \leq \frac{|I_{<n}|}{|I_n \setminus S|} \leq \frac{|I_{<n}|}{|I_n|/2} < 2^{-n+1} \leq \frac{1}{2} - \varepsilon.$$

As $n$ can be arbitrarily large, $p$ cannot force that $S$ $\varepsilon$-almost bisects $\dot{X}$. $\square$

Applying this lemma and the aforementioned properties of $\mathbb{P}_{\omega_2}$, we obtain the following:

**Theorem 4.2.** *If $V \vDash \text{CH}$, then $V^{\mathbb{P}_{\omega_2}} \vDash \omega_1 = \text{cof}(\mathcal{M}) < \mathfrak{s}_{1/2\pm\varepsilon} = \mathfrak{s}_{1/2} = \omega_2$.*

We show that this strict inequality can be obtained via a finite support iteration as well, namely, we can dualise (now in the forcing sense) the result saying that $\mathfrak{r}_{1/2\pm\varepsilon} < \text{add}(\mathcal{M})$ in the Hechler model (see Theorem 3.4). Consider the dual Hechler model, that is, every model of the form $V^{\mathbb{D}_{\omega_1}}$ where $V \vDash$ "$\text{MA} + \mathfrak{c} \geq \omega_2$". We know (see e. g. [BJ95, Model 7.6.10]) that $\text{cof}(\mathcal{M}) = \omega_1$ holds in these models. We interpret $\mathbb{D}$ as the filter-based Laver forcing for the Fréchet filter over $\omega^{\uparrow<\omega} = \{s \in \omega^{<\omega} \mid s \text{ is strictly increasing}\}$, that is, $T \in \mathbb{D}$ if $T \subseteq \omega^{\uparrow<\omega}$ is a tree (i. e. $T$ is closed with regard to taking initial segments) with a fixed element $\text{stem}(T) \in T$ such that

(a) either $t \subseteq \text{stem}(T)$ or $\text{stem}(T) \subseteq t$ for every $t \in T$ and
(b) $\text{ext}_T(t) = \{n \in \omega \mid t^\frown(n) \in T\}$ is cofinite in $\omega$ for every $t \in T$ with $\text{stem}(T) \subseteq t$.

If $T \in \mathbb{D}$ and $s \in T$, then let $T{\restriction_s} = \{t \in T \mid t \subseteq s \text{ or } t \supseteq s\} \in \mathbb{D}$. Let $\dot{d}$ be a $\mathbb{D}$-name for the generic dominating real, i. e. $\dot{d} = \bigcup\{\text{stem}(T) \mid T \text{ belongs to the generic filter}\} \in \omega^{\uparrow\omega} = \{f \in \omega^\omega \mid f \text{ is strictly increasing}\}$.

We recall a classical preservation theorem we will apply (see e. g. [Gol93] or [BJ95, section 6.4]). Fix a sequence $(\sqsubset_n)_{n\in\omega}$ of increasing closed relations on $\omega^\omega$ such that

$$(\sqsubset_n)^g = \{f \in \omega^\omega \mid f \sqsubset_n g\} \text{ is nowhere dense}$$

for every $n$ and $g$. Let $\sqsubset = \bigcup_{n\in\omega} \sqsubset_n$. We will use terminology compatible with the one we use when working with relational systems: If $\kappa = \text{cof}(\kappa) > \omega$, then a $U \subseteq \omega^\omega$ is $\kappa$-$\sqsubset$-*unbounded* if for every $C \subseteq \omega^\omega$ of size $< \kappa$, there is an $f \in U$ which is not $\sqsubset$-below any element of $C$. In this case we will write $f \not\sqsubset C$.

**Definition 4.3.** Let $\kappa = \text{cof}(\kappa) > \omega$ and let $\mathbb{P}$ be a $\kappa$-cc forcing notion. We say that $\mathbb{P}$ is $\kappa$-$\sqsubset$-*good* if for every $\mathbb{P}$-name $\dot{h}$ for an element of $\omega^\omega$, there exists a non-empty $Y \subseteq \omega^\omega$ of size $<\kappa$ such that $\Vdash_{\mathbb{P}} f \not\sqsubset \dot{h}$ whenever $f \not\sqsubset Y$. Say that $\mathbb{P}$ is $\sqsubset$-*good* if it is $\omega_1$-$\sqsubset$-good.[1]

It is straightforward to show that if $\mathbb{P}$ is $\kappa$-$\sqsubset$-good, then $\mathbb{P}$ preserves (a) "$F$ is $\kappa$-$\sqsubset$-unbounded" for $F \subseteq \omega^\omega$ and (b) "$\mathfrak{d}(\sqsubset) \geq \lambda$" for cardinals $\lambda \geq \kappa$ (see [Mej13]).

**Theorem 4.4** (see [JS90] or [BJ95, Lemma 6.4.12])**.** *Let $\kappa = \text{cof}(\kappa) > \omega$ and assume that $(\mathbb{P}_\alpha, \dot{\mathbb{Q}}_\beta)_{\alpha\leq\delta,\beta<\delta}$ is a finite support iteration of $\kappa$-cc forcing notions such that $\Vdash_\alpha$ "$\dot{\mathbb{Q}}_\alpha$ is $\kappa$-$\sqsubset$-good" for every $\alpha < \delta$. Then $\mathbb{P}_\delta$ is $\kappa$-$\sqsubset$-good as well.*

We are going to apply this theorem with $\dot{\mathbb{Q}}_\alpha = \mathbb{D}$, but first we define our relation $\sqsubset$. Fix an $\varepsilon \in (0, 1/2)$ and an interval partition $(I_k)$ as we already did above, that is, $|I_0| \geq 2$ and $|I_k| > 2^k |I_{<k}|$ for every $k$. Let

$$\mathcal{O} = \{(H, (E_k)_{k\in\omega}) \mid H \in [\omega]^\omega, E_k \subseteq I_k \text{ and } |E_k|/|I_k| > 1/4 \text{ for every } k\}$$

and define the relation $\sqsubset = \bigcup_{n\in\omega} \sqsubset_n$ on $[\omega]^\omega \times \mathcal{O}$ by

$$X \sqsubset_n (H, (E_k)_{k\in\omega}) \iff \forall\, k \in H \setminus n \colon |X \cap E_k| < (1/2 + \varepsilon)\big(|X \cap I_k| + |I_{<k}|\big).$$

First of all, notice that $[\omega]^\omega$ as a subspace of $\mathcal{P}(\omega)$ is canonically homeomorphic (denoted by $\simeq$) to $\omega^{\uparrow\omega} \subseteq \omega^\omega$ and hence to $\omega^\omega$ itself. To code $\mathcal{O}$ as $\omega^\omega$ as well, let

$$\mathcal{Q}_k = \left\{ E \subseteq I_k \;\middle|\; \frac{|E|}{|I_k|} > \frac{1}{4} \right\};$$

[1] See [Mej13] for the proof that this variant of $\kappa$-$\sqsubset$-goodness is equivalent (under the assumption that $\mathbb{P}$ is $\kappa$-cc) to the classical definition from [JS90] or [BJ95, Definition 6.4.4].

if we consider $\mathcal{Q}_k$ as a discrete space, then $\prod_{k\in\omega}\mathcal{Q}_k$ is a compact metric space and

$$\mathcal{O} = [\omega]^\omega \times \prod_{k\in\omega}\mathcal{Q}_k \simeq \omega^\omega \times \prod_{k\in\omega}\mathcal{Q}_k \simeq \prod_{k\in\omega}(\omega\times\mathcal{Q}_k) \simeq \omega^\omega.$$

**Lemma 4.5.** *The following statements hold:*

(1) $\operatorname{dom}(\sqsubset) = [\omega]^\omega$ *and* $\operatorname{ran}(\sqsubset) = \mathcal{O}$.

(2) *The relation* $\sqsubset_n$ *is closed in* $[\omega]^\omega \times \mathcal{O}$.

(3) *The set* $(\sqsubset_n)^{(H,(E_k))}$ *is nowhere dense in* $[\omega]^\omega$ *for every* $(H,(E_k)) \in \mathcal{O}$.

(4) $\operatorname{cov}(\mathcal{M}) \leq \mathfrak{d}(\sqsubset) \leq \mathfrak{s}_{1/2\pm\varepsilon}$.

*Proof.* (1): First let $X \in [\omega]^\omega$. If there is an infinite $H \subseteq \omega$ such that $|X\cap I_k|/|I_k| < 3/4$ for every $k \in H$, then $X \sqsubset_1 (H, (I_k \setminus X)_{k\in\omega})$. If there is no such $H$, then $|X \cap I_k|/|I_k| \geq 3/4$ for every $k \geq K$ for some $K \in \omega$, and if we choose some $E_k \subseteq I_k$ with $|E_k|/|I_k| \in (1/4, 3/8)$ for every $k \geq K (\geq 2)$, then

$$|X \cap E_k| \leq |E_k| < \frac{3}{8}|I_k| < \Big(\frac{1}{2}+\varepsilon\Big)|X\cap I_k| \text{ for every } k \geq K,$$

hence $X \sqsubset_0 (\omega\setminus K, (E_k)_{k\in\omega})$.

Now let $(H, (E_k)_{k\in\omega}) \in \mathcal{O}$. If $X = \omega \setminus \bigcup_{k\in H} E_k$ is infinite, then $X \sqsubset_1 (H, (E_k))$. If this set is finite, however, then there is a $K$ such that $\omega\setminus K \subseteq H$ and $E_k = I_k$ for every $k \geq K$. If $Y \in [\omega]^\omega$ such that $|Y \cap I_k| = 1$ for every $k$, then

$$|Y \cap E_k| \leq 1 < \Big(\frac{1}{2}+\varepsilon\Big)\big(1 + |I_{<k}|\big) \text{ holds for every } k \in H\setminus 1,$$

and so $Y \sqsubset_1 (H, (E_k))$.

(2): If $X \not\sqsubset_n (H, (E_k))$, then this is witnessed by a $k_0 \in H \setminus n$, that is,

$$|X \cap E_{k_0}| \geq \Big(\frac{1}{2}+\varepsilon\Big)\big(|X\cap I_{k_0}| + |I_{<k_0}|\big).$$

Then $X' \not\sqsubset_n (H', (E'_k))$ for all $X' \in [\omega]^\omega$ and $(H', (E'_k)) \in \mathcal{O}$ such that $X' \cap I_{k_0} = X \cap I_{k_0}$, $k_0 \in H'$ and $E'_{k_0} = E_{k_0}$; these pairs $(X', (H', (E'_k))) \in [\omega]^\omega \times \mathcal{O}$ form an open neighbourhood of $(X, (H, (E_k)))$ in $[\omega]^\omega \times \mathcal{O}$.

(3): Fix an $X \in [\omega]^\omega$ with $X \sqsubset_n (H, (E_k))$ and a basic open neighbourhood $U_m = \{Y \in [\omega]^\omega \mid Y \cap m = X \cap m\}$ of $X$ (for an $m \in \omega$). If $k \in H \setminus n$ such that $\min(I_k) \geq m$ (and $k$ is sufficiently large, see below), and $Y \in U_m$ such that $Y \cap I_k = E_k$, then $k$ witnesses $Y \not\sqsubset_n (H, (E_k))$, that is, $|Y \cap E_k| = |E_k| > (1/2 + \varepsilon)(|E_k| + |I_{<k}|)$ because

$$\frac{|E_k|}{|I_k|} > \Big(\frac{1}{2}+\varepsilon\Big)\Big(\frac{|E_k|}{|I_k|} + 2^{-k}\Big)$$

holds if $k$ is sufficiently large.

(4): $\operatorname{cov}(\mathcal{M}) \leq \mathfrak{d}(\sqsubset)$ follows from (1) and (3). Now let $S \in [\omega]^\omega$, define $S' = S$ if $|S \cap I_k|/|I_k| > 1/4$ for infinitely many $k$ and $S' = \omega \setminus S$ otherwise, and let $H_S = \{k \mid |S' \cap I_k|/|I_k| > 1/4\} \in [\omega]^\omega$. For $k \in H_S$, let $E^S_k = S' \cap I_k$, and for $k \in \omega \setminus H_S$, let $E_k = I_k$. It is enough to show that if $S$ $\varepsilon$-almost bisects $X$, then $X \sqsubset (H_S, (E^S_k))$, i.e. that $\big(\mathrm{id}, S \longmapsto (H_S, (E^S_k))\big)\colon ([\omega]^\omega, \sqsubset, \mathcal{O}) \to \mathbf{Reap}^\perp_{1/2\pm\varepsilon}$ is a Tukey connection. Clearly, $S \mid_{1/2\pm\varepsilon} X$ iff $S' \mid_{1/2\pm\varepsilon} X$, therefore if $k \in H_S$ is sufficiently large, then $X \cap I_{\leq k} \neq \varnothing$ and

$$\frac{|X \cap E^S_k|}{|X\cap I_k| + |I_{<k}|} \leq \frac{|S' \cap X \cap I_{\leq k}|}{|X \cap I_{\leq k}|} < \frac{1}{2}+\varepsilon,$$

finishing the proof. $\square$

**Lemma 4.6.** $\mathbb{D}$ *is* $\sqsubset$*-good.*

*Proof.* Let $(\dot{H}, (\dot{E}_k))$ be a $\mathbb{D}$-name for an element of $\mathcal{O}$. We will construct a countable family $\mathcal{O}' \subseteq \mathcal{O}$ such that whenever $X \in [\omega]^\omega \cap V$ and $X \not\sqsubset \mathcal{O}'$, then $\Vdash X \not\sqsubset (\dot{H}, (\dot{E}_k))$.

Let $\dot{H} = \{\dot{k}_0 < \dot{k}_1 < \dots\}$ be an enumeration in $V^{\mathbb{D}}$. Recall that we denote the generic real of $\mathbb{D}$ by $\dot{d}$. By thinning out $\dot{H}$, we can assume that $\dot{d}(n) < \dot{k}_n$ for every $n$. (Note that if $\Vdash X \not\sqsubset (\dot{J}, (\dot{E}_k))$ for some infinite $\dot{J} \subseteq \dot{H}$, then this holds for $\dot{H}$ as well.)

We define a rank function $\rho_n$ on $\omega^{<\uparrow\omega}$ for every fixed $n \in \omega$ as follows: We set $\rho_n(s) = 0$ if there are $k_{n,s} \in \omega$ and $E_{n,s} \subseteq I_{k_{n,s}}$ such that whenever $T \in \mathbb{D}$ and $\operatorname{stem}(T) = s$, then there is a $T' \leq T$ which forces that "$\dot{k}_n = k_{n,s}$ and $\dot{E}_{k_{n,s}} = E_{n,s}$". Then we proceed by recursion: At the $\alpha$th stage, after already having defined $\{s \in \omega^{<\omega} \mid \rho_n(s) = \beta\}$ for every $\beta < \alpha$, we set $\rho_n(s) = \alpha$ if

$$Y_{n,s} = \{m \mid \rho_n\big(s^\frown(m)\big) < \alpha\} \text{ is infinite.}$$

We show that $\operatorname{dom}(\rho_n) = \omega^{<\uparrow\omega}$ for every $n$. Assume towards a contradiction that $\rho_n(s)$ is not defined. Then $\{m \mid \rho_n(s^\frown(m)) \text{ is not defined}\}$ is cofinite, hence we can construct a $T \in \mathbb{D}$ with stem $s$ such that $\rho_n(t)$ is not defined for every $t \in T$ above $s$. There are a $T' \leq T$, $k \in \omega$ and $E \subseteq I_k$ such that $T' \Vdash$ "$\dot{k}_n = k$ and $\dot{E}_k = E$"; in particular, $k$ and $E$ witness that $\rho_n(\operatorname{stem}(T')) = 0$, a contradiction.

Also, we will need the fact that if $n \geq |s|$, then $\rho_n(s) > 0$. To see this, let (for $k \in \omega$) $T_k \in \mathbb{D}$ such that $\operatorname{stem}(T_k) = s$ and $\operatorname{ext}_{T_k}(s) = \omega \setminus k$. Then $T_k \Vdash k \leq \dot{d}(|s|) \leq \dot{d}(n) < \dot{k}_n$. If $\rho_n(s) = 0$, then we could pick $k_{n,s}$ such that $T \Vdash \dot{k}_n = k_{n,s}$ whenever $T \in \mathbb{D}$ with stem $s$; but in that case, $T = T_{k_{n,s}} \Vdash k_{n,s} < \dot{k}_n$, a contradiction.

Now, if $\rho_n(s) = 1$ and $m \in Y_{n,s}$, then $\rho_n(s^\frown(m)) = 0$ and hence we have defined $k_{n,s^\frown(m)}$ and $E_{n,s^\frown(m)}$. Note that

$$\{m \in Y_{n,s} \mid k_{n,s^\frown(m)} = k\} \text{ is finite}$$

for each $k$. Otherwise, there are $k \in \omega$ and $E \subseteq I_k$ such that $X = \{m \in Y_{n,s} \mid k_{n,s^\frown(m)} = k \text{ and } E_{n,s^\frown(m)} = E\}$ is also infinite, and hence if $T \in \mathbb{D}$ with $\operatorname{stem}(T) = s$, then there is an $m \in X \cap \operatorname{ext}_T(s)$, and so there is a $T' \leq T\restriction_{s^\frown(m)} \leq T$ such that $T' \Vdash$ "$\dot{k}_n = k$ and $\dot{E}_k = E$", in other words, $\rho_n(s) = 0$, a contradiction.

For such $n$ and $s$ we can thus fix an infinite $Z_{n,s} \subseteq Y_{n,s}$ such that

$$\text{if } m, m' \in Z_{n,s} \text{ and } m < m', \text{ then } k_{n,s^\frown(m)} < k_{n,s^\frown(m')}.$$

We let $K_{n,s} = \{k_{n,s^\frown(m)} \mid m \in Z_{n,s}\}$, $E^{n,s}_k = E_{n,s^\frown(m)}$ if $m \in Z_{n,s}$ and $k_{n,s^\frown(m)} = k$ and $E^{n,s}_k = I_k$ otherwise, and

$$\mathcal{O}' = \{(K_{n,s}, (E^{n,s}_k)_{k\in\omega}) \mid n \in \omega,\ s \in \omega^{\uparrow<\omega},\ \rho_n(s) = 1\} \subseteq \mathcal{O}.$$

To finish the proof, fix an $X \in [\omega]^\omega \cap V$ and assume that $X \not\sqsubset \mathcal{O}'$, i. e. $X \not\sqsubset (K_{n,s}, (E^{n,s}_k))$ whenever $\rho_n(s) = 1$, or more explicitly, for infinitely many $k \in K_{n,s}$

$$(\bullet_{k,X,(E^{n,s}_k)}) \qquad |X \cap E^{n,s}_k| \geq \Big(\frac{1}{2} + \varepsilon\Big)\big(|X \cap I_k| + |I_{<k}|\big).$$

To show that $\Vdash X \not\sqsubset (\dot{H}, (\dot{E}_k))$, we fix $T \in \mathbb{D}$ and $n \in \omega$. We will find a $T' \leq T$ and a $k \geq n$ such that $T' \Vdash$ "$k \in \dot{H}$ and $(\bullet_{k,X,(\dot{E}_k)})$". We can assume that $n \geq |\operatorname{stem}(T)|$ and hence $\rho_n(\operatorname{stem}(T)) > 0$. By induction on this rank, one can easily show that there is an $s \in T$ above the stem such that $\rho_n(s) = 1$. Pick a $k = k_{n,s^\frown(m)} \in K_{n,s} \cap \operatorname{ext}_T(s) \setminus n$ such that $(\bullet_{k,X,(E^{n,s}_k)})$. By the definition of $\rho_n(s^\frown(m)) = 0$, there is a $T' \leq T\restriction_{s^\frown(m)} \leq T$ which forces that $\dot{k}_n = k$ and $\dot{E}_k = E_{n,s^\frown(m)} = E^{n,s}_k$; in particular, $T'$ also forces that $k \in \dot{H} \setminus n$ and $(\bullet_{k,X,(\dot{E}_k)})$. $\square$

Applying Lemma 4.6 and Theorem 4.4, we obtain that $\mathbb{D}_\delta$ preserves $\mathfrak{s}_{1/2\pm\varepsilon} \geq \mathfrak{d}(\sqsubset) \geq \operatorname{cov}(\mathcal{M})^V$, and hence the following:

**Theorem 4.7.** *If $V \vDash \mathrm{MA} + \mathfrak{c} \geq \omega_2$, then $V^{\mathbb{D}_{\omega_1}} \vDash \omega_1 = \operatorname{cof}(\mathcal{M}) < \mathfrak{s}_{1/2\pm\varepsilon} = \mathfrak{s}_{1/2} = \mathfrak{c}$.*

## 5. Further Generalisations: $\mathfrak{s}_\rho$ and $\mathfrak{r}_\rho$

In this last section, we take a look at generalisations of $\mathfrak{s}_{1/2}$ and $\mathfrak{r}_{1/2}$, in the following sense: For $\rho \in (0,1)$, let $\mathbf{Reap}_\rho = ([\omega]^\omega, \not|_\rho, [\omega]^\omega)$ where $S \mid_\rho X$ ("$S$ *$\rho$-splits* $X$") if

$$\frac{|S \cap X \cap n|}{|X \cap n|} \xrightarrow{n \to \infty} \rho.$$

We write $\mathfrak{s}_\rho = \mathfrak{b}(\mathbf{Reap}_\rho)$ and $\mathfrak{r}_\rho = \mathfrak{d}(\mathbf{Reap}_\rho)$.

Obviously, $\mathbf{Reap}_{1-\rho} \equiv \mathbf{Reap}_\rho \preccurlyeq \mathbf{Reap}$ and hence $\mathfrak{s}_\rho = \mathfrak{s}_{1-\rho} \geq \mathfrak{s}$ and $\mathfrak{r}_\rho = \mathfrak{r}_{1-\rho} \leq \mathfrak{r}$ for all $\rho \in (0,1)$. We will need the following easy observation:

**Fact 5.1.** Let $\rho_0, \rho_1 \in (0,1)$ and $A, B, X \in [\omega]^\omega$. Then the following statements hold:

(1) If $A \mid_{\rho_0} X$ and $B \mid_{\rho_1} A \cap X$, then $A \cap B \mid_{\rho_0\rho_1} X$; hence $\max\{\mathfrak{s}_{\rho_0}, \mathfrak{s}_{\rho_1}\} \geq \mathfrak{s}_{\rho_0\rho_1}$ and $\min\{\mathfrak{r}_{\rho_0}, \mathfrak{r}_{\rho_1}\} \leq \mathfrak{r}_{\rho_0\rho_1}$.
(2) If $A \mid_{\rho_0} X$ and $B \mid_{\rho_1} X \setminus A$, then $A \cup B \mid_{\rho_0+\rho_1-\rho_0\rho_1} X$; hence $\max\{\mathfrak{s}_{\rho_0}, \mathfrak{s}_{\rho_1}\} \geq \mathfrak{s}_{\rho_0+\rho_1-\rho_0\rho_1}$ and $\min\{\mathfrak{r}_{\rho_0}, \mathfrak{r}_{\rho_1}\} \leq \mathfrak{r}_{\rho_0+\rho_1-\rho_0\rho_1}$.

*Proof.* (1) follows from

$$\frac{|(A \cap B) \cap X \cap n|}{|X \cap n|} = \frac{|B \cap (A \cap X) \cap n|}{|(A \cap X) \cap n|} \cdot \frac{|A \cap X \cap n|}{|X \cap n|},$$

and (2) follows from

$$\begin{aligned}\frac{|(A \cup B) \cap X \cap n|}{|X \cap n|} &= \frac{|A \cap X \cap n| + |B \cap (X \setminus A) \cap n|}{|X \cap n|} \\ &= \frac{|A \cap X \cap n|}{|X \cap n|} + \frac{|B \cap (X \setminus A) \cap n|}{|(X \setminus A) \cap n|} \cdot \frac{|(X \setminus A) \cap n|}{|X \cap n|}.\end{aligned}$$ □

Also, we recall a classical construction:

**Fact 5.2** (non-integer bases)**.** Let $b > 1$ be a real number. Then every $x > 0$ can be written as $x = \sum_{n=-\infty}^{N} c_n b^n$ where $N \geq 0$ is an integer and $0 \leq c_n < b$ are also integers for all $n$.

**Theorem 5.3.** *$\mathfrak{r}_\rho$ does not depend on $\rho$.*

*Proof.* Fix a $\rho \in (1/2, 1)$. We show that $\mathfrak{r}_\rho = \mathfrak{r}_{1/2}$.

To show that $\mathfrak{r}_\rho \geq \mathfrak{r}_{1/2}$, let $\mathcal{R} \subseteq [\omega]^\omega$ such that $|\mathcal{R}| < \mathfrak{r}_{1/2}$. We will construct an $S \in [\omega]^\omega$ such that $S \mid_\rho \mathcal{R}$ (that is, such that $S \mid_\rho X$ for every $X \in \mathcal{R}$).

By recursion on $n \in \omega$, one can easily construct $S_n \in [\omega]^\omega$ and $\mathcal{R}_n \subseteq [\omega]^\omega$ such that $S_0 = \omega$, $\mathcal{R}_0 = \mathcal{R}$, $S_{n+1} \mid_{1/2} \mathcal{R}_n$ and

$$\mathcal{R}_{n+1} = \{S_{n+1} \cap X \mid X \in \mathcal{R}_n\}.$$

For $m \geq 1$, define

$$I_m = \bigcap_{n \leq m} S_n \qquad \text{and} \qquad D_m = I_{m-1} \setminus I_m = I_{m-1} \setminus S_m.$$

It follows that $\mathcal{R}_m = \{I_m \cap X \mid X \in \mathcal{R}\}$ and hence that $I_m, D_m \in [\omega]^\omega$ for every $m \geq 1$. First of all, we show that $I_m \mid_{1/2^m} \mathcal{R}$ and $D_m \mid_{1/2^m} \mathcal{R}$ for every $m \geq 1$. This holds for $m = 1$ by the definitions above. We proceed by induction on $m$; assume that the claim holds for a fixed $m$. If $X \in \mathcal{R}$, then $I_m \cap X \in \mathcal{R}_m$ and thus $S_{m+1} \mid_{1/2} I_m \cap X$. Since $I_m \mid_{1/2^m} X$ holds as well, we can apply Fact 5.1 (1) with $\rho_0 = 1/2^m$, $\rho_1 = 1/2$, $A = I_m$, $B = S_{m+1}$ and $X$, and obtain that $I_{m+1} = S_{m+1} \cap I_m \mid_{1/2^{m+1}} X$. It follows that $D_{m+1} = I_m \setminus S_{m+1} \mid_{1/2^{m+1}} X$ (because $1/2^m - 1/2^{m+1} = 1/2^{m+1}$).

Now let $P \subseteq \omega \setminus \{0\}$ be such that $\sum_{m\in P} 2^{-m} = \rho$ (a representation of $\rho < 1$ in base 2) and $S = \bigcup_{m\in P} D_m$. We show that $S \mid_\rho \mathcal{R}$. Fix an $X \in \mathcal{R}$ and define the *lower* and *upper relative density of* $S \in [\omega]^\omega$ *in* $X \in [\omega]^\omega$ as follows:

$$\overline{d}_X(S) = \limsup_{n\to\infty} \frac{|S\cap X\cap n|}{|X\cap n|},$$

$$\underline{d}_X(S) = \liminf_{n\to\infty} \frac{|S\cap X\cap n|}{|X\cap n|}.$$

Clearly, $\overline{d}_X$ and $\underline{d}_X$ are monotone, $\overline{d}_X(\omega\setminus S) = 1-\underline{d}_X(S)$ and $\underline{d}_X(\omega\setminus S) = 1-\overline{d}_X(S)$, and $S \mid_\rho X$ iff $\overline{d}_X(S) = \underline{d}_X(S) = \rho$; in this case we write $d_X(S) = \rho$. Also, notice that $d_X$ is finitely additive. It follows that it is enough to show that $\underline{d}_X(S) \geq \rho$ and $\underline{d}_X(\omega\setminus S) \geq 1-\rho$. For $k\in\omega$, let $A_k = \bigcup_{m\in P\cap k} D_m$ (and $\bigcup\varnothing = \varnothing$), then $\underline{d}_X(A_k) \leq \underline{d}_X(S)$, and at the same time $\underline{d}_X(A_k) = d_X(A_k) = \sum_{m\in P\cap k} d_X(D_m) = \sum_{m\in P\cap k} 2^{-m} \xrightarrow{k\to\infty} \rho$, hence $\rho \leq \underline{d}_X(S)$. To show that $\underline{d}_X(\omega\setminus S) \geq 1-\rho$, write $\omega\setminus S$ as $\bigcup_{m\in\omega\setminus P} D_m$. (Since finite modifications of $S_n$ do not affect the above, we can assume without loss of generality that $\bigcap_{n\in\omega} S_n = \varnothing$ and hence that $\omega = \bigcup_{m\geq 1} D_m$ is a partition.)

The proof of the converse inequality is similar but with a little twist: We start with a family $\mathcal{R}' \subseteq [\omega]^\omega$ such that $|\mathcal{R}'| < \mathfrak{r}_\rho$ and have to construct an $S' \in [\omega]^\omega$ such that $S' \mid_{1/2} \mathcal{R}'$. By recursion on $n\in\omega$, we define $S'_n \in [\omega]^\omega$ and $\mathcal{R}'_n \subseteq [\omega]^\omega$ such that $S'_0 = \omega$, $\mathcal{R}'_0 = \mathcal{R}'$, $S'_{n+1} \mid_\rho \mathcal{R}'_n$ and $\mathcal{R}'_{n+1} = \{S'_{n+1}\cap X \mid X \in \mathcal{R}'_n\}$. For $m\geq 1$, define $I'_m = \bigcap_{n\leq m} S_n$ and $D'_m = I_{m-1}\setminus I_m$. It follows (by induction) that $I'_m \mid_{\rho^m} \mathcal{R}'$ and $D'_m \mid_{\rho^{m-1}(1-\rho)} \mathcal{R}'$ for every $m \geq 1$. If we can find a $P' \subseteq \omega\setminus\{0\}$ such that $\sum_{m\in P'} \rho^{m-1}(1-\rho) = 1/2$, then, just like above, $S' = \bigcup_{m\in P'} D'_m \mid_{1/2} \mathcal{R}'$. The problem is that when writing ${}^{\rho}/_{2-2\rho}$ in the non-integer base $1 < b = \rho^{-1} < 2$, i. e. ${}^{\rho}/_{2-2\rho} = \sum_{n=-\infty}^{N} (\rho^{-1})^n$, we may have to use non-zero coefficients for at least some $0 \leq n \leq N$. Therefore, we have two cases:

Case 1: If ${}^{\rho}/_{2-2\rho} < 1$ i. e. $\rho < 2/3$, then $N < 0$ and we can pick an appropriate $P'$ and hence conclude that $\mathfrak{r}_\rho \leq \mathfrak{r}_{1/2}$.

Case 2: If $2/3 \leq \rho < 1$ then, applying Fact 5.1 (1), $\mathfrak{r}_\rho = \min\{\mathfrak{r}_\rho, \mathfrak{r}_\rho\} \leq \mathfrak{r}_{\rho^2} = \min\{\mathfrak{r}_{\rho^2}, \mathfrak{r}_{\rho^2}\} \leq \mathfrak{r}_{\rho^4} \leq \dots$; in other words, $\mathfrak{r}_\rho \leq \mathfrak{r}_{\rho^2} \leq \mathfrak{r}_{\rho^4} \leq \mathfrak{r}_{\rho^8} \leq \dots$ There is an $n$ such that $1/3 < \rho^{2^n} < 2/3$ (because $1/3 < (2/3)^2$), and hence $\mathfrak{r}_\rho \leq \mathfrak{r}_{\rho^{2^n}} = \mathfrak{r}_{1-\rho^{2^n}} \leq \mathfrak{r}_{1/2}$. □

Note that we did not make use of Tukey connections in the proof above. Therefore, the dual equality of $\rho$-splitting numbers for different parameters $\rho$ does not follow. In fact, we only showed that a family too small to be $1/2$-reaping is not $\rho$-reaping and vice versa – which does not, in any obvious way, yield a method to turn a $1/2$-reaping family into a $\rho$-reaping family or vice versa.

**Question D.** *Is it consistent that* $\mathfrak{s}_\rho \neq \mathfrak{s}_\tau$ *for some* $\rho, \tau \in (0,1)$*?*

The next natural question is if we can generalise $\mathbf{Cov}(\mathcal{N}) \preccurlyeq \mathbf{Reap}_{1/2} \preccurlyeq \mathbf{Cov}(\mathcal{M})^\perp$, that is, the inequalities $\operatorname{non}(\mathcal{N}) \geq \mathfrak{s}_{1/2} \geq \operatorname{cov}(\mathcal{M})$ and $\operatorname{cov}(\mathcal{N}) \leq \mathfrak{r}_{1/2} \leq \operatorname{non}(\mathcal{M})$ for arbitrary $\rho \in (0,1)$.

$\mathbf{Cov}(\mathcal{N}) \preccurlyeq \mathbf{Reap}_\rho$ is problematic because in the case of $\rho = 1/2$, the proof uses the law of large numbers to show that $\{S \in [\omega]^\omega \mid S \mid_{1/2} R\}$ is of measure 1 for every fixed $R \in [\omega]^\omega$.

**Question E.** *Does* $\mathbf{Cov}(\mathcal{N}) \preccurlyeq \mathbf{Reap}_\rho$, *or at least* $\mathrm{non}(\mathcal{N}) \geq \mathfrak{s}_\rho$ *and* $\mathrm{cov}(\mathcal{N}) \leq \mathfrak{r}_\rho$, *hold?*

However, $\mathbf{Reap}_\rho \preccurlyeq \mathbf{Cov}(\mathcal{M})^\perp$ and hence $\mathfrak{s}_\rho \geq \mathrm{cov}(\mathcal{M})$ and $\mathfrak{r}_\rho \leq \mathrm{non}(\mathcal{M})$ hold because it is easy to see that $\{X \in [\omega]^\omega \mid S \mid_\rho X\} =$

$$\bigcap_{\varepsilon>0}\bigcup_{N\in\omega}\bigcap_{n\geq N}\{X \in [\omega]^\omega \mid (\rho-\varepsilon)|X\cap n| \leq |S\cap X\cap n| \leq (\rho+\varepsilon)|X\cap n|\} \in \mathcal{M}$$

and hence, identifying $\mathcal{P}(\omega)$ and $2^\omega$, $(F,G)\colon ([\omega]^\omega, \nmid_\rho, [\omega]^\omega) \to (\mathcal{M}, \not\ni, \mathcal{P}(\omega))$ with $F(S) = [\omega]^{<\omega} \cup \{X \in [\omega]^\omega \mid S \mid_\rho X\}$ and $G(X) = X$ if $X$ is infinite and $G(X) = \omega$ if $X$ is finite is a (Borel) Tukey connection.

We can define $\mathfrak{s}_0$, $\mathfrak{r}_0$, $\mathfrak{s}_1$, and $\mathfrak{r}_1$ as well. To avoid the trivial case $\mathfrak{s}_1 = 1$ and to maintain the duality $\mathfrak{s}_1 = \mathfrak{s}_0$ and $\mathfrak{r}_1 = \mathfrak{r}_0$, we say that $S \mid_0 R$ ("$S$ 0-*splits* $R$") if $S$ is infinite and coinfinite, $R$ is infinite and $|S\cap R\cap n|/|R\cap n| \to 0$; we define $S \mid_1 R$ similarly. Hence let $\mathbf{Reap}_0 = (\{S \in [\omega]^\omega \mid |\omega \setminus S| = \omega\}, \mid_0, [\omega]^\omega)$, $\mathbf{Reap}_1 = (\{S \in [\omega]^\omega \mid |\omega \setminus S| = \omega\}, \mid_1, [\omega]^\omega)$, $\mathfrak{s}_0 = \mathfrak{b}(\mathbf{Reap}_0) = \mathfrak{s}_1 = \mathfrak{b}(\mathbf{Reap}_1)$ and $\mathfrak{r}_0 = \mathfrak{d}(\mathbf{Reap}_0) = \mathfrak{r}_1 = \mathfrak{d}(\mathbf{Reap}_1)$.

Just like for $\rho \in (0,1)$, if $S \in [\omega]^\omega$ and $|\omega \setminus S| = \omega$, then $\{X \in [\omega]^\omega \mid S \mid_0 X\} \in \mathcal{M}$ and hence $\mathbf{Reap}_0 \preccurlyeq \mathbf{Cov}(\mathcal{M})^\perp$; in particular, $\mathfrak{s}_0 \geq \mathrm{cov}(\mathcal{M})$ and $\mathfrak{r}_0 \leq \mathrm{non}(\mathcal{M})$.

**Fact 5.4.** $\mathbf{Dom}^\perp \preccurlyeq \mathbf{Reap}_0$ and hence $\mathfrak{d} \geq \mathfrak{s}_0$ and $\mathfrak{b} \leq \mathfrak{r}_0$.

*Proof.* Instead of $\omega^\omega$, we work with $\mathcal{X} = \{x \in \omega^{\uparrow\omega} \mid |\omega \setminus \mathrm{ran}(x)| = \omega\}$. It is trivial to show that $(\mathcal{X}, \leq^*, \mathcal{X}) \equiv (\omega^\omega, \leq^*, \omega^\omega) = \mathbf{Dom}$. We define a Tukey connection $(F,G)\colon \mathbf{Reap}_0^\perp \to \mathbf{Dom}$, that is, an

$$(F,G)\colon \big([\omega]^\omega, \text{"is 0-split by"}, \{S \in [\omega]^\omega \mid |\omega \setminus S| = \omega\}\big) \to (\mathcal{X}, \leq^*, \mathcal{X})$$

as follows: Let $F\colon [\omega]^\omega \to \mathcal{X}$ be defined by $F(R)(n) = r_{2^n}$ where $R = \{r_0 < r_1 < \dots\} \in [\omega]^\omega$ and let $G\colon \mathcal{X} \to \{S \in [\omega]^\omega \mid |\omega \setminus S| = \omega\}$ be defined by $G(x) = \mathrm{ran}(x)$. If $F(R) \leq^* x$, $r_{2^n} \leq x(n)$ for every $n \geq N$ for some $N$ and $r_{2^n} < k \leq r_{2^{n+1}}$, then

$$\frac{|\mathrm{ran}(x)\cap R\cap k|}{|R\cap k|} \leq \frac{N+n}{2^n} \xrightarrow{n\to\infty} 0,$$

hence $G(x) = \mathrm{ran}(x)$ 0-splits $R$. □

**Question F.** *Does* $\mathbf{Dom}^\perp \equiv \mathbf{Reap}_0$, *and hence* $\mathfrak{d} = \mathfrak{s}_0$ *and* $\mathfrak{b} = \mathfrak{r}_0$, *hold? Or do at least* $\mathfrak{s}_0 \geq \mathfrak{s}$ *and* $\mathfrak{r}_0 \leq \mathfrak{r}$ *hold?*

Institute of Discrete Mathematics and Geometry, TU Wien, Wiedner Hauptstrasse 8–10/104, 1040 Vienna, Austria

*Email address*: barnabasfarkas@gmail.com

Institute of IT Security Research and Center for Artificial Intelligence, St. Pölten University of Applied Sciences, Campus-Platz 1, 3100 St. Pölten, Austria

*Email address*: mail@l17r.eu

*URL*: https://l17r.eu

Institute of Mathematics, University of Zurich, Winterthurerstrasse 190, 8057 Zurich, Switzerland

*Email address*: marc.lischka@math.uzh.ch, marc.lischka@gmail.com